\documentclass[fleqn,onecolumn,leqno]{amsproc}%
\usepackage{amsfonts}
\usepackage{amsmath}
\usepackage{amssymb}
\usepackage{graphicx}
\usepackage[none]{hyphenat}%
\theoremstyle{plain}

\numberwithin{equation}{section}
\begin{document}
\title{}

\begin{center}
\bigskip{\LARGE \textbf{Ratio of Price to Expectation and}}

{\LARGE \textbf{Complete Bernstein Functions}}

\bigskip

\textbf{Yukio Hirashita}

\bigskip

\textbf{Abstract\smallskip}
\end{center}

\noindent For a game with positive expectation and some negative profit, a
unique price exists, at which the optimal proportion of investment reaches its
maximum. For a game with parallel translated profit, the ratio of this price
to its expectation tends to converge toward less than or equal to $1/2$ if its
expectation converges to $0^{+}$. In this paper, we will investigate such
properties by using the integral representations of a complete Bernstein
function and establish several Abelian and Tauberian
theorems.\textbf{\smallskip\smallskip}



\bigskip

\begin{center}
\textbf{1. Introduction\smallskip}
\end{center}

\noindent Consider a coin-flipping game such that profit is $9$ dollars or
$-2$ dollars if a tossed coin yields heads or tails, respectively. For
simplicity, we will omit the currency notation. Let $t\in\lbrack0,$ $1]$ be
the proportion of investment. Then, the investor repeatedly invests $t$ of
his/her current capital (see [12, 13]). For example, let $c>0$ be the current
capital; when the investor plays the game once, his/her capital will be
$9ct/u+c(1-t)$ or $-2ct/u+c(1-t)$ if a tossed coin yields heads or tails,
respectively, where $u>0$ is the price of the game such that $u/(u+2)>t$. Let
the initial capital be $1$. After $N$ attempts, if the investor has capital
$c_{N}$, then the growth rate (geometric mean) is given by $c_{N}^{1/N}$. As
the value $G_{u}(t):=\lim_{N\rightarrow\infty}\left(  \text{expectation of
}c_{N}^{1/N}\right)  =$ $\sqrt{(9t/u+1-t)(-2t/u+1-t)}$ is a function with
respect to $t$, it reaches its maximum at $t=t_{u}$ $=(7/2-u)u$
$/((u+2)(9-u))$. It is noteworthy that the value $\lim_{N\rightarrow\infty
}\left(  \text{variance of }c_{N}^{1/N}\right)  $ is $0$.

In general, a game $(a(x),$ $F(x))$ would mean that if the investor invests
$1$ unit (which price is $u$ dollars), then he/she receives $a(x)$ dollars
(including the invested money) in accordance with a distribution function
$F(x),$ defined on an interval $I\subseteq(-\infty$, $\infty)$ such that
$\int_{I}d(F(x))=1$. It is assumed that the profit function $a(x)$ is
measurable and non-constant (a.e.) with respect to $F(x)$. When no confusion
arises, we write $dF$ for $d(F(x))$ and use the following notation:%
\begin{equation}
E:=\int_{I}a(x)dF,\text{ }\xi:=\mathrm{ess}\text{ }\inf_{x\in I}\,a(x),\text{
}H_{\xi}:=\int_{I}\frac{1}{a(x)-\xi}dF. \tag{1.1}%
\end{equation}

In this paper, we always assume that $E>0$ and $\xi>-\infty$. If
$\int_{a(x)=\xi}dF>0$, we define $H_{\xi}$ $=\infty$ and $1/H_{\xi}$ $=0.$
Since $a(x)$ is non-constant, we have $\xi<E$, $H_{\xi}>0$, $1/H_{\xi}<\infty
$, and $\xi+1/H_{\xi}<E$.

In order to explain the background of this paper, we will define notations
such as $w_{\beta}(z)$ and $G_{u}(t)$ in this paragraph. However, this paper
utilizes neither such notations nor their related properties, except in the
first paragraph of Section 2. We denote the integral $\int_{I}(a(x)-\beta
)/(a(x)z-z\beta+\beta)dF$ by $w_{\beta}(z)$, which is holomorphic with respect
to two complex variables $(z,$ $\beta)$ ($z:=t+si,$ $\beta:=u+hi$,
$i:=\sqrt{-1}$, $\{t,s,u,h\}\subset%
\mathbb{R}
$) near each point ($t_{0},$ $u_{0}$) such that $0<t_{0}<u_{0}/(u_{0}-\xi)$
and $u_{0}>\max(0,$ $\xi)$. We denote $\exp(\int_{I}\log\left(
a(x)t/u-t+1\right)  dF)$ by $G_{u}(t)$ and term it as \textit{the limit
expectation of growth rate }for each $u>0$ and $0\leq t\leq1$ with $\xi
t/u-t+1>0$. We say that $t_{u}$ is \textit{the optimal proportion of
investment} with respect to $u>0$, if
\begin{equation}
\overline{\lim}_{\substack{\rho\rightarrow t_{u}\\0\leq\rho\leq1\\\xi
\rho/u-\rho+1>0}}\int_{I}\log\frac{a(x)t/u-t+1}{a(x)\rho/u-\rho+1}dF\leq0
\tag{1.2}%
\end{equation}
for each $0\leq t\leq1$ with $\xi t/u-t+1>0.$ A game $(a(x)$, $F(x))$ is said
to be \textit{effective} if $\int_{a(x)>1}a(x)^{\nu}dF<\infty$ for some
$\nu>0$. If a game is effective, $G_{u}(t)$ is continuous (see [9, Theorem
4.1]) and the inequality (1.2) implies that $G_{u}(t_{u})=\sup_{0\leq
t\leq1,\text{ }\xi t/u-t+1>0}G_{u}(t)$, which suggests that $t_{u}$ is optimal
for maximizing the limit expectation of growth rate.

For a game with\textit{ parallel translated profit} $(a(x)-m,$ $F(x))$
$(m<E)$, we use underlined notations such as $\underline{a}(x):=a(x)-m,$
$\underline{E}:=E-m,$ $\underline{\xi}:=\xi-m,$ and $\underline{H}%
_{\underline{\xi}}:=H_{\xi}$.

From [9, Lemma 3.16], if $m\in(\xi,$ $E)$, then a unique price $\underline
{u}_{\max}\in(0,$ $E-m)$ exists such that $\underline{t}_{u}$ is strictly
increasing in the interval $0<u<\underline{u}_{\max}$ and strictly decreasing
in the interval $\underline{u}_{\max}<u<E-m.$ It should be noted that
$\underline{u}_{\max}$ is a function with respect to $m\in(\xi,$ $E),$ and it
satisfies $\underline{t}_{\underline{u}_{\max}}=\max_{0<u<E-m}\underline
{t}_{u}$. In a sense, $\underline{u}_{\max}$ is considered to be the price in
which the broker's commission income is maximized.

Under mild restrictions, we will show that $\lim_{m\rightarrow E^{-}%
}\underline{u}_{\max}/\underline{E}=1/2$ (see Theorem 3.19). In such a case,
it suggests that the so-called half price sale makes a profit. For example, in
the case of the abovementioned coin-flipping game, we obtain%
\[
\lim_{m\rightarrow E^{-}}\frac{\underline{u}_{\max}}{\underline{E}}%
=\lim_{m\rightarrow(7/2)^{-}}\frac{11\sqrt{(m+2)(9-m)}/2-(m+2)(9-m)}%
{(7/2-m)^{2}}=\frac{1}{2},
\]
where $-2<m<7/2$ (see Corollary 3.20).

Defining $\Psi(c):=1/\int_{I}(a(x)+c)^{-1}dF-c$ $(c\in(-\xi,$ $\infty))$, we
obtain the following:

L{\small EMMA} 1.1.\textbf{ }$\lim_{c\rightarrow\infty}\Psi(c)=E.$

P{\small ROOF}.\textbf{ }Assume $c>\max(1,-2\xi)$. Then, we have $a(x)+c/2>0$
and $0<\int_{I}c/(a(x)$ $+c)dF$ $\leq\int_{I}2dF=2.$ If $E<\infty$, then, by
applying Lebesgue's monotone convergence and dominated convergence
theorems\ to the equation%
\[
\frac{1}{\int_{I}\frac{1}{a(x)+c}dF}-c=E+\frac{E^{2}}{c-E}-\frac{\int_{I}%
\frac{a(x)^{2}}{a(x)+c}dF}{\left(  1-\frac{E}{c}\right)  \int_{I}\frac
{c}{a(x)+c}dF}\text{ \ \ }(c\neq E),
\]
we obtain the conclusion (even if $E\leq0$). Assume $E=\infty$. Since $a(x)$
is non-constant with respect to $F(x)$, we observe that%
\[
\Psi^{\prime}(c)=\frac{\int_{I}\frac{1}{(a(x)+c)^{2}}dF-\left(  \int_{I}%
\frac{1}{a(x)+c}dF\right)  ^{2}}{\left(  \int_{I}\frac{1}{a(x)+c}dF\right)
^{2}}>0,
\]
which implies that $\Psi(c)$ is increasing with respect to $c$. Putting
$\lim_{c\rightarrow\infty}\Psi(c)=M$ (including $\infty$),%
\[
a_{N}(x):=\left\{
\begin{array}
[c]{cc}%
N, & a(x)>N,\\
a(x), & a(x)\leq N,
\end{array}
\right.  \text{ and\ }b_{c,N}:=\frac{1}{\int_{I}\frac{1}{a_{N}(x)+c}%
dF}-c\text{ \ \ \ }(N>\max(1,\text{ }\xi)).
\]
Then, the following properties hold:

(1) $a_{N}(x)$ is nondecreasing with respect to $N$. (2) $b_{c,N}$ is
nondecreasing with respect to $N$. (3) From the above arguments, we obtain
$\lim_{c\rightarrow\infty}b_{c,N}$ $=\int_{I}a_{N}(x)dF$ $<\infty$, which is
nondecreasing with respect to $N$. (4) By applying Lebesgue's (monotone
convergence) theorem, we obtain $\lim_{N\rightarrow\infty}b_{c,N}=\Psi(c)$,
which is increasing with respect to $c$. (5) Further, from Lebesgue's
(monotone convergence) theorem,we obtain $\lim_{N\rightarrow\infty}%
(\lim_{c\rightarrow\infty}b_{c,N})=\lim_{N\rightarrow\infty}\int_{I}%
a_{N}(x)dF$ $=\int_{I}a(x)dF=\infty.$

Therefore, if $M<\infty$, then $\Psi(c)\leq M$ and $b_{c,N}\leq M$, which
contradicts the fact that $\lim_{N\rightarrow\infty}(\lim_{c\rightarrow\infty
}b_{c,N})=\infty$. Hence, $M=E=\infty$, that is, $\lim_{c\rightarrow\infty
}\Psi(c)=E$.

\hfill$\square$

\bigskip

\begin{center}
\textbf{2. Parallel translated profit\smallskip}
\end{center}

We consider a game with parallel translated profit $(a(x)-m,$ $F(x))$ to have
\textit{sufficiently small positive expectation}, if $\xi+1/H_{\xi}<m<E$. In
this case, it is easy to observe that $\underline{E}=E-m>0$, $\underline{\xi
}=\xi-m<0$ and $\underline{\xi}+1/\underline{H}_{\underline{\xi}}=\xi-m$
$+1/H_{\xi}<0$. Therefore, from [9, Lemma 4.27], $\eta_{m}:=\lim
_{u\rightarrow0^{+}}\underline{t}_{u}/u$ exists such that $0<\eta
_{m}<-1/\underline{\xi}=1/(m-\xi)$. For each $u\in(0,$ $E-m)$ and $t\in(0,$
$u/(u-\underline{\xi}))$, we have $\underline{w}_{u}(\underline{t}_{u})=0$
where $\underline{w}_{u}(t)=$ $\int_{I}\left(  a(x)-m-u\right)  $ $/\left(
(a(x)-m)t-tu+u\right)  dF$ and $\underline{t}_{u}\in(0,u/(u-\underline{\xi}%
))$. It should be noted that $(a(x)-m)\underline{t}_{u}/u-\underline{t}%
_{u}+1>0$ and $(a(x)-m)\eta_{m}+1$ $\geq1-(m-\xi)\eta_{m}$ $>0$ for each $x\in
I$. The equation $\underline{w}_{u}(\underline{t}_{u})=0$ can be written as
$\int_{I}\left(  (a(x)-m)\underline{t}_{u}/u-\underline{t}_{u}+1\right)
^{-1}dF$ $=1$. Hence, we have%
\begin{equation}
\int_{I}\frac{1}{(a(x)-m)\eta_{m}+1}dF=1, \tag{2.1}%
\end{equation}
because on the set $\{x$ $|$ $a(x)\geq m\}$, $1/((a(x)-m)\underline{t}%
_{u}/u-\underline{t}_{u}+1)$ is strictly increasing with respect to
sufficiently small $u>0$ (see [9, Lemmas 3.12, 3.15 and 3.16]), and on the set
$\{x$ $|$ $\xi\leq a(x)<m\},$ $1/((a(x)-m)\underline{t}_{u}/u-\underline
{t}_{u}+1)$ converges uniformly to $1/((a(x)-m)\eta_{m}+1)$\ $(u\rightarrow
0^{+})$.

Since $\Psi(c)$ $(-\xi<c<\infty)$ is strictly increasing from $\xi+1/H_{\xi}$
to $E$ (see the proof of Lemma 1.1), the equation $m=\Psi(c)$ has a unique
solution $c=c_{m}$ for each $\xi+1/H_{\xi}<m<E$. Since the equation $m$
$=1/\int_{I}(a(x)+c_{m})^{-1}dF$ $-c_{m}$ is equivalent to%
\[
\int_{I}\frac{1}{(a(x)-m)\frac{1}{m+c_{m}}+1}dF=1,
\]
from (2.1), we obtain that $\eta_{m}=1/(m+c_{m})$.

L{\small EMMA} 2.1. $c_{m}$ \textit{is strictly increasing from} $-\xi$
\textit{to} $\infty$ \textit{with respect to} $m\in(\xi+1/H_{\xi}$, $E).$

P{\small ROOF}. As $\Psi(c)$ is strictly increasing from $\xi+1/H_{\xi}$ to
$E$, the relation $m=\Psi(c_{m})$ leads to the conclusion.\hfill$\square$

L{\small EMMA} 2.2. $\eta_{m}$ \textit{is strictly decreasing from} $H_{\xi}$
\textit{to} $0$ \textit{with respect to} $m\in(\xi+1/H_{\xi}$, $E).$

P{\small ROOF}. Since $\eta_{m}=1/(m+c_{m}),$ Lemma 2.1 leads to the
conclusion.\hfill$\square$

L{\small EMMA} 2.3. $\lim_{m\rightarrow E^{-}}m\eta_{m}=0$\textbf{.}

P{\small ROOF}. From $\lim_{m\rightarrow E^{-}}c_{m}=\infty$ and Lebesgue's
theorem, we obtain the equality $m\eta_{m}=m/(m+c_{m})=\int_{I}%
a(x)/(a(x)+c_{m})dF$, which implies the conclusion.

\hfill$\square$

L{\small EMMA} 2.4. $\underline{t}_{u}=u\eta_{m+u}$ $(0<u<E-m).$

P{\small ROOF}. From the property of $\eta_{m+u}$, we observe that $\int
_{I}\left(  (a(x)-(m+u))\eta_{m+u}+1\right)  ^{-1}$ $dF=1$, which can be
written as $\int_{I}\left(  (a(x)-m)(u\eta_{m+u})/u-u\eta_{m+u}+1\right)
^{-1}dF=1$. This suggests that $\underline{w}_{u}(u\eta_{m+u})=0$. Therefore,
by the uniqueness of $\underline{t}_{u}$, we arrive at the conclusion.\hfill
$\square$

L{\small EMMA} 2.5. $\underline{u}_{\max}$ \textit{can be uniquely determined
by the system}%
\begin{equation}
\left\{
\begin{array}
[c]{c}%
m+v=\Psi(c),\\
v=(m+c)\Psi^{\prime}(c),
\end{array}
\right.  \tag{2.2}%
\end{equation}
\textit{with two unknown variables }$v$\textit{ }$(=\underline{u}_{\max}%
)$\textit{ and }$c$\textit{ }$(=c_{m+\underline{u}_{\max}})$\textit{, for each
}$m\in(\xi+1/H_{\xi},$\textit{ }$E)$\textit{.}

P{\small ROOF}. From $m=\Psi(c_{m})$, we obtain $m+\underline{u}_{\max}%
=\Psi(c_{m+\underline{u}_{\max}})$. Since $\underline{t}_{\underline{u}_{\max
}}$ $=\max_{0<u<E-m}\underline{t}_{u}$, we find that $\underline
{t}_{\underline{u}_{\max}}^{\prime}=0.$ From $\underline{t}_{u}=u\eta
_{m+u}=u/(m+u+c_{m+u})$, we obtain $\underline{t}_{u}^{\prime}=\left(
m+c_{m+u}-uc_{m+u}^{\prime}\right)  /(m+u+c_{m+u})^{2}$. Thus,
$m+c_{m+\underline{u}_{\max}}$ $-\underline{u}_{\max}c_{m+\underline{u}_{\max
}}^{\prime}=0.$ Using $\Psi^{\prime}(c_{m})c_{m}^{\prime}=1$, we have
$m+c_{m+\underline{u}_{\max}}$ $-\underline{u}_{\max}/\Psi^{\prime
}(c_{m+\underline{u}_{\max}})$ $=0$, which implies (2.2). On the other hand,
from%
\[
\Psi^{\prime\prime}(c)=\frac{2\left(  \left(  \int_{I}\frac{1}{(a(x)+c)^{2}%
}dF\right)  ^{2}-\int_{I}\frac{1}{a(x)+c}dF\int_{I}\frac{1}{(a(x)+c)3}%
dF\right)  }{\left(  \int_{I}\frac{1}{a(x)+c}dF\right)  ^{3}}<0\text{
\ \ \ }(c>-\xi),
\]
$\Psi(c)$\ is a strictly concave function due to Schwarz's inequality.
Therefore, a line $y-m=\Psi^{\prime}(c)(x+m)$ that is tangent to $\Psi(c)$ and
passes through the point $(-m$, $m)$ is uniquely determined. This implies the
uniqueness of the solution of (2.2).\hfill$\square$

E{\small XAMPLE} 2.6. The game $(x,$ $\int_{0}^{x}1/(\pi(t+1)\sqrt{t})dt)$
$(x\in(0,\infty))$ has the following properties:\textit{ }$\xi=0,$
$\xi+1/H_{\xi}=0<m<E=\infty,$ $\Psi(c)=\sqrt{c}$, $\eta_{m}=1/(m(m+1))$,
$c_{m}=m^{2}$, $\underline{t}_{u}$ $=u/((m+u)(m+u+1))$, and $\underline
{u}_{\max}=\sqrt{m(m+1)}$.

E{\small XAMPLE} 2.7. The game $(x,$ $\int_{0}^{x}8r\sqrt{rt}/(\pi
(t+r)^{3})dt)$ $(x\in(0,\infty),$ $r>0)$ has the following properties:
$\xi=0,$ $\xi+1/H_{\xi}=r/3<m<E=3r,$ $\Psi(c)$ $=(c-r)^{3}$ $/(c^{2}%
-6cr-3r^{2}+8r\sqrt{cr})$ $-c$, $\eta_{m}=(3r-m)^{2}/(m+r)^{3}$, $c_{m}$
$=r(r-3m)^{2}/(3r-m)^{2}$, $\underline{t}_{u}$ $=u(3r-m-u)^{2}$ $/(r+m+u)^{3}%
$, and $\underline{u}_{\max}=(3r-m)(m+r)$ $/(m+9r)$.

R{\small EMARK}.\textbf{ }Since $\underline{E}=E-m>0$, the assumption $E>0$
can be dismissed as long as we consider a game $(a(x)-m,$ $F(x))$ with
$\xi+1/H_{\xi}<m<E$. It is clear that even if $E\leq0$, $\xi+1/H_{\xi}<E$
holds, provided $\xi=\mathrm{ess}$ $\inf_{x\in I}\,a(x)>-\infty$.

\newpage

\begin{center}
3\textbf{. Complete Bernstein functions\smallskip}
\end{center}

A $C^{\infty}$ function $f$ $:$ $(0,$ $\infty)\rightarrow%
\mathbb{R}
$ with a continuous extension to $[0$, $\infty)$ is called a \textit{Bernstein
function} if $f\geq0$ and $(-1)^{k}f^{\,(k)}(x)\leq0$ for each $k\in
\{1,2,3,...\}$ (see [5, Definition 1.2.1]).

A function $f:$ $(0$, $\infty)\rightarrow%
\mathbb{R}
$ is called a\textit{ complete Bernstein function} if there exists a Bernstein
function\textit{ }$\phi$ such that $f(x)=x^{2}\int_{0}^{\infty}e^{-sx}%
\phi(s)ds$ (see [14, Definition 1.4]).

T{\small HEOREM} 3.1 [14, Theorem 1.5]. \textit{Each of the following five
properties of }$f$\textit{ }$:$\textit{ }$(0,$\textit{ }$\infty)\rightarrow%
\mathbb{R}
$\textit{ implies the other four}:

(1) $f$ \textit{is a complete Bernstein function}.

(2) $f$ \textit{can be represented as }$f(x)=\tau x+b+\int_{0}^{\infty
}x/(x+t)\sigma(dt)$\textit{ with }$\tau,b\geq0$

\ \ \ \ \textit{and a measure }$\sigma$\textit{ on }$(0,$ $\infty).$

(3) $f$ \textit{extends analytically on }$%
\mathbb{C}
\backslash(-\infty,$\textit{ }$0]$\textit{ such that }$f(\overline
{z})=\overline{f(z)}$\textit{ and }$\operatorname{Im}z\operatorname{Im}f(z)$

\ \ \ \ $\geq0.$ (\textit{In other words,} $f$\textit{ preserves the upper and
lower half-planes in }$%
\mathbb{C}
$).

(4) $f$ \textit{is a Bernstein function with representation }$f(x)$\textit{
}$=\tau x+b$\textit{ }

\ \ \ \ \textit{ }$+\int_{0}^{\infty}(1-e^{-sx})\beta(s)ds$, \textit{where
}$\tau,b\geq0,$\textit{ }$\beta(s)=\int_{0}^{\infty}e^{-st}\rho(dt),$\textit{
and }

\ \ \ \ \textit{ }$\int_{0}^{\infty}1/\left(  t(t+1)\right)  \rho(dt)<\infty.$
(\textit{In fact, }$\rho(dt)=t\sigma(dt)$\textit{ of} (2).)

(5) $x/f(x)$ \textit{is a complete Bernstein function or }$f\equiv0$.

Note that the triple $(\tau,$ $b,$ $\rho)$ given above is uniquely determined
by $f$ (see [5, Theorem 1.2.3]).

L{\small EMMA} 3.2\textbf{. }\textit{The function }$\Psi(c)=1/\int
_{I}(a(x)+c)^{-1}dF-c$: $(-\xi,$\textit{ }$\infty)\rightarrow(\xi$ $+1/H_{\xi
},$\textit{ }$E)$\textit{ extends analytically on }$%
\mathbb{C}
\backslash(-\infty,-\xi]$\textit{ and preserves the upper and lower
half-planes.}

P{\small ROOF}. From Lemma 1.1, we obtain $\Psi((-\xi,$ $\infty))\subset
(\xi+1/H_{\xi},$ $E)$. Putting $c$ $=u+yi\in%
\mathbb{C}
\backslash(-\infty,-\xi]$, $\{u,$ $y\}\subset%
\mathbb{R}
$ and $s=a(x)+u,$ then
\[
\Psi(c)=\frac{1}{\int_{I}\frac{1}{s+yi}dF}-u-yi=\frac{1}{\int_{I}\frac
{s}{s^{2}+y^{2}}dF-yi\int_{I}\frac{1}{s^{2}+y^{2}}dF}-u-yi.
\]
If $y\gtrless0$, then due to Schwarz's inequality, we observe that%
\[
\operatorname{Im}\Psi(c)=\frac{y\left(  \int_{I}\frac{1}{s^{2}+y^{2}}%
dF\int_{I}\frac{s^{2}}{s^{2}+y^{2}}dF-\left(  \int_{I}\frac{1}{\sqrt
{s^{2}+y^{2}}}\frac{s}{\sqrt{s^{2}+y^{2}}}dF\right)  ^{2}\right)  }{\left(
\int_{I}\frac{s}{s^{2}+y^{2}}dF\right)  ^{2}+y^{2}\left(  \int_{I}\frac
{1}{s^{2}+y^{2}}dF\right)  ^{2}}\gtrless0.
\]

Set $\alpha(v):=\int_{a(x)\leq v+\xi}dF.$ Then, $\alpha(v)$ is a right
continuous nondecreasing function such that $\alpha(v)\geq0,$ $\alpha(\xi
^{-})=0$ and $\alpha(\infty)=1$. Thus, the Stieltjes transform $\int
_{I}(a(x)+c)^{-1}dF=\int_{0^{-}}^{\infty}(v+c+\xi)^{-1}d(\alpha(v))$ is
analytic with respect to $t=c+\xi\in%
\mathbb{C}
\backslash(-\infty,0]$ (see [15, Corollary VIII.2b.1]). It is easy to verify
that $\Psi(c)$ has no singular point in $%
\mathbb{C}
\backslash(-\infty,-\xi]$.\hfill$\square$

T{\small HEOREM} 3.3.\textbf{ }$\Psi(c)-\xi-1/H_{\xi}$\textit{ is a complete
Bernstein function with respect to }$t=c+\xi>0$.

P{\small ROOF}. From Theorem 3.1 (3) and Lemma 3.2, we arrive at the conclusion.

\hfill$\square$

L{\small EMMA} 3.4. $\lim_{c\rightarrow\infty}\Psi(c)/c=0$.

P{\small ROOF}. From $\Psi(c)/c=1/(1-\int_{I}a(x)/(a(x)+c)dF)-1$ and%
\[
\lim_{c\rightarrow\infty}\int_{I}\frac{a(x)}{a(x)+c}dF=\lim_{c\rightarrow
\infty}\int_{a(x)>0}\frac{a(x)}{a(x)+c}dF+\lim_{c\rightarrow\infty}\int
_{\xi\leq a(x)<0}\frac{a(x)}{a(x)+c}dF=0,
\]
we arrive at the conclusion by applying Lebesgue's monotone convergence and
dominated convergence theorems.\hfill$\square$

L{\small EMMA} 3.5. $\Psi(c)$ \textit{can be written as}%
\begin{equation}
\Psi(c)=\xi+\frac{1}{H_{\xi}}+\int_{0}^{\infty}\frac{c+\xi}{t(t+c+\xi)}%
\rho(dt)\text{ \ \ \ (}c>-\xi\text{)} \tag{3.1}%
\end{equation}
\textit{with }$\int_{0}^{\infty}t^{-1}\rho(dt)=E-\xi-1/H_{\xi}$\textit{ and
}$\int_{0}^{\infty}1/\left(  t(t+1)\right)  \rho(dt)<\infty$.

P{\small ROOF}. From Theorems 3.1 and 3.3,\textbf{ }$\Psi(c)$ can be written
as
\[
\Psi(c)-\xi-1/H_{\xi}=\tau(c+\xi)+b+\int_{0}^{\infty}\frac{c+\xi}{t(t+c+\xi
)}\rho(dt)\text{ \ \ \ (}c>-\xi\text{),}%
\]
where $\tau\geq0$, $b\geq0,$ and $\int_{0}^{\infty}1/\left(  t(t+1)\right)
\rho(dt)<\infty$. Since
\[
\frac{\Psi(c)}{c}=\tau+\frac{\tau\xi+b+\xi+1/H_{\xi}}{c}+\frac{c+\xi}{c}%
\int_{0}^{\infty}\frac{1}{t(t+c+\xi)}\rho(dt),
\]
we have $\tau=0$ by applying Lemma 3.4 and Lebesgue's theorem. Since
$\partial\left(  (c+\xi)/(t(t+c+\xi))\right)  $ $/\partial c$ $=(t+c+\xi
)^{-2}$ $>0$, $(c+\xi)/(t(t+c+\xi))$ is increasing with respect to $c>-\xi$.
From $\Psi(-\xi)$ $=\xi+1/H_{\xi}$ and Lebesgue's theorem, we obtain $b=0.$
From $\lim_{c\rightarrow\infty}\Psi(c)=E$ and Lebesgue's theorem, we obtain
that $E-\xi-1/H_{\xi}$ $=\int_{0}^{\infty}t^{-1}\rho(dt)$.$\hfill\square$

L{\small EMMA} 3.6. \textit{The condition }$E<\infty$\textit{ is equivalent to
}$\int_{0}^{\infty}t^{-1}\rho(dt)<\infty$.

P{\small ROOF}. Lemma 3.5 shows that $\int_{0}^{\infty}t^{-1}\rho
(dt)=E-\xi-1/H_{\xi}$, which implies the conclusion because $\xi$ and
$1/H_{\xi}$ are finite.$\hfill\square$

L{\small EMMA} 3.7. \textit{A function }$f(x)\geq0$\textit{ }$(x>0)$\textit{
is a complete Bernstein function if and only if }$1/(x+f(x))$\textit{ is a
Stieltjes transform. In this case, a right continuous nondecreasing function
}$0\leq G(t)\leq1$\textit{\ exists such that} $1/(x+f(x))$ $=\int_{0^{-}%
}^{\infty}(x$ $+t)^{-1}d(G(t))$.

P{\small ROOF}. Let $f(x)$ be a complete Bernstein function. Therefore, in
accordance with Theorem 3.1 (2), $x+f(x)$ is a complete Bernstein function.
Further, in accordance with Theorem 3.1 (5), $x/(x+f(x))$ is a complete
Bernstein function. Thus, from Theorem 3.1 (2), we have $x/(x+f(x))=\widetilde
{\tau}x+\widetilde{b}+\int_{0}^{\infty}x/(x+t)\widetilde{\sigma}(dt)$ $(x>0)$
with $\widetilde{\tau},$ $\widetilde{b}\geq0$ and a measure $\widetilde
{\sigma}\ $on ($0$, $\infty$). We can obtain $\widetilde{\tau}=0$ as follows.
From Theorem 3.1 (2), we have $f(x)=\tau x+b+\int_{0}^{\infty}x/(x+t)\sigma
(dt)$ $(\tau,$ $b\geq0)$. Thus,
\[
\frac{1}{\widetilde{\tau}x+\widetilde{b}+\int_{0}^{\infty}\frac{x}%
{x+t}\widetilde{\sigma}(dt)}=1+\tau+\frac{b}{x}+\int_{0}^{\infty}\frac{1}%
{x+t}\sigma(dt).
\]
If $\widetilde{\tau}>0,$ then the process $x\rightarrow\infty$ leads to
$0=1+\tau$, which contradicts the fact that $\tau\geq0.$ Therefore, the right
continuous nondecreasing function $G(t):=\widetilde{b}+\int_{0}^{t}%
\widetilde{\sigma}(dt)$ (if $x\geq0$) or 0 (if $x<0$) yields the Stieltjes
transform%
\[
\frac{1}{x+f(x)}=\frac{\widetilde{b}}{x}+\int_{0}^{\infty}\frac{1}%
{x+t}\widetilde{\sigma}(dt)=\int_{0^{-}}^{\infty}\frac{1}{x+t}d(G(t)).
\]
If $\int_{0^{-}}^{\infty}d(G(t))>1,$ then $\lim\inf_{x\rightarrow\infty}f(x)$
$=\lim\inf_{x\rightarrow\infty}$ $\ x(1/\int_{0^{-}}^{\infty}x(x+t)^{-1}%
d(G(t))$ $-1)$ $<0$, which contradicts the assumption that $f(x)\geq0$. Thus,
we find that $0\leq G(t)\leq1$.

On the other hand, assume that $1/(x+f(x))$ is a Stieltjes transform such that
$1/(x+f(x))=\widehat{\tau}+\frac{\widehat{b}}{x}+\int_{0}^{\infty}%
(x+t)^{-1}\widehat{\sigma}(dt)$, where $\widehat{\tau}$ and $\widehat{b}$\ are
constants, and $\widehat{\sigma}$\ is a measure on $(0$, $\infty)$. Since
$f(x)\geq0$, by applying Lebesgue's theorem, we obtain $\lim_{x\rightarrow
\infty}1/(x+f(x))=0=\widehat{\tau}.$ Put $G(t)=\widehat{b}+\int_{0}%
^{t}\widehat{\sigma}(dt)$ $($if $x\geq0)$ or $0$ $($if $x<0)$. Then, we obtain
$f(x)=1/\int_{0^{-}}^{\infty}(x+t)^{-1}d(G(t))-x.$ As mentioned above,
$\int_{0^{-}}^{\infty}d(G(t))>1$ causes a contradiction. Thus, we have
$\int_{0^{-}}^{\infty}(dG(t))\leq1$. Putting $x=u+yi$ and $s=t+u$, as in the
proof of Lemma 3.2, we observe that if $y\gtrless0$, then%
\begin{align}
&  \operatorname{Im}f(u+yi)\tag{3.2}\\
&  =\frac{y\left(
\begin{array}
[c]{c}%
\left(  1-\int_{0^{-}}^{\infty}dG\right)  \int_{0^{-}}^{\infty}\frac{1}%
{s^{2}+y^{2}}dG+\int_{0^{-}}^{\infty}\frac{1}{s^{2}+y^{2}}dG\int_{0^{-}%
}^{\infty}\frac{s^{2}}{s^{2}+y^{2}}dG\\
-\left(  \int_{0^{-}}^{\infty}\frac{1}{\sqrt{s^{2}+y^{2}}}\frac{s}{\sqrt
{s^{2}+y^{2}}}dG\right)  ^{2}%
\end{array}
\right)  }{\left(  \int_{0^{-}}^{\infty}\frac{s}{s^{2}+y^{2}}dG\right)
^{2}+y^{2}\left(  \int_{0^{-}}^{\infty}\frac{1}{s^{2}+y^{2}}dG\right)  ^{2}%
}\nonumber\\
&  \gtrless0.\nonumber
\end{align}
Thus, the analytic function $f(u+yi)$ on $%
\mathbb{C}
\backslash(-\infty,$ $0]$ (see [15, Corollary VIII.2b.1]) preserves the upper
and lower half-planes. This implies, in accordance with Theorem 3.1 (3), that
$f(x)$ is a complete Bernstein function.$\hfill\square$

We characterize the relation between the subset of complete Bernstein
functions such that $\tau=0$ and all the probability measures on $[0,$
$\infty)$.

T{\small HEOREM} 3.8. \textit{A complete Bernstein function }$f(x)=\tau
x+b+\int_{0}^{\infty}x/(x+t)\sigma(dt)$\textit{ can be written as
}$f(x)=1/\left(  \int_{0^{-}}^{\infty}(x+t)^{-1}d(G(t))\right)  -x$\textit{
with a distribution function }$0\leq G(t)\leq1$\textit{ with }$G(\infty
)=1$\textit{ and }$G(0^{-})=0$\textit{ if and only if }$\tau=0$\textit{.}

P{\small ROOF}. From Lemma 3.7, for a complete Bernstein function $\tau x+b$
$+\int_{0}^{\infty}x/(x+t)\sigma(dt)$, a right continuous nondecreasing
function $0\leq G(t)\leq1$ exists such that $1/\left(  x+\tau x+b+\int
_{0}^{\infty}x/(x+t)\sigma(dt)\right)  $ $=\int_{0^{-}}^{\infty}%
(x+t)^{-1}d(G(t))$. Therefore, we observe that $\int_{0^{-}}^{\infty
}x/(x+t)d(G(t))=1/\left(  1+\tau+b/x+\int_{0}^{\infty}(x+t)^{-1}%
\sigma(dt)\right)  ,$ which implies $\int_{0^{-}}^{\infty}d(G(t))=1/(1+\tau)$
as $x\rightarrow\infty.$ Thus, if $\tau=0$, then $\int_{0^{-}}^{\infty
}d(G(t))=1$. In this case, we have $f(x)$ $=1/\left(  \int_{0^{-}}^{\infty
}(x+t)^{-1}d(G(t))\right)  $ $-x$. The converse is obtained by applying
Lebesgue's theorem to the equation $1/\int_{0^{-}}^{\infty}%
x/(x+t)d(G(t))-1=\tau+b/x$ $\ +\int_{0}^{\infty}(x+t)^{-1}\sigma(dt).$%
\hfill$\square$

L{\small EMMA} 3.9. $\lim_{c\rightarrow\infty}\Psi^{(n)}(c)=0$ $(n=1,2,3,...)$.

P{\small ROOF}. Using (3.1), we have (see [15, Corollary VIII.2b.2])%
\begin{align}
\Psi^{\prime}(c)  &  =\int_{0}^{\infty}\frac{1}{(t+c+\xi)^{2}}\rho
(dt),\Psi^{\prime\prime}(c)=-\int_{0}^{\infty}\frac{2}{(t+c+\xi)^{3}}%
\rho(dt),...\tag{3.3}\\
\Psi^{(n)}(c)  &  =(-1)^{n-1}n!\int_{0}^{\infty}\frac{1}{(t+c+\xi)^{n+1}}%
\rho(dt).\nonumber
\end{align}
By applying Lebesgue's (monotone convergence) theorem, we conclude that
$\lim_{c\rightarrow\infty}\Psi^{(n)}(c)$ $=0$ $(n=1,2,3,...)$.\hfill$\square$

L{\small EMMA} 3.10. \textit{If }$E<\infty$\textit{, then }$\lim
_{c\rightarrow\infty}c^{n}\Psi^{\text{ }(n)}(c)=0$\textit{ }$(n=1,2,3,...)$%
\textit{.}

P{\small ROOF}. From Lemma 3.6, $\int_{0}^{\infty}t^{-1}\rho(dt)<\infty$. From
properties such as
\begin{align*}
c^{n}\Psi^{(n)}(c)  &  =(-1)^{n-1}n!\frac{c^{n}}{(c+\xi)^{n}}\int_{0}^{\infty
}\frac{t(c+\xi)^{n}}{(t+c+\xi)^{n+1}}\frac{\rho(dt)}{t},\\
\lim_{c\rightarrow\infty}\frac{c^{n}}{(c+\xi)^{n}}  &  =1,\text{ }%
\lim_{c\rightarrow\infty}\frac{t(c+\xi)^{n}}{(t+c+\xi)^{n+1}}=0\text{,}\\
\left\vert \frac{t(c+\xi)^{n}}{(t+c+\xi)^{n+1}}\right\vert  &  <1\text{
\ \ \ \ }(t,\text{ }c+\xi>0),
\end{align*}
we can apply Lebesgue's (dominated convergence) theorem and obtain
$\lim_{c\rightarrow\infty}c^{n}\Psi^{(n)}(c)$ $=0$ $(n=1,2,3,...)$%
.\hfill$\square$

L{\small EMMA} 3.11. $\lim_{m\rightarrow E^{-}}\underline{u}_{\max}%
=\lim_{c\rightarrow\infty}c\Psi^{\prime}(c)$\textit{ if one of them exists. In
particular, if }$E<\infty$, $\lim_{m\rightarrow E^{-}}\underline{u}_{\max}=0$.

P{\small ROOF}. From (2.2), we obtain
\begin{equation}
\underline{u}_{\max}=\frac{(1+\frac{\Psi(c_{m+\underline{u}_{\max}}%
)}{c_{m+\underline{u}_{\max}}})}{1+\Psi^{\prime}(c_{m+\underline{u}_{\max}}%
)}c_{m+\underline{u}_{\max}}\Psi^{\prime}(c_{m+\underline{u}_{\max}}).
\tag{3.4}%
\end{equation}
From $\lim\inf_{m\rightarrow E^{-}}c_{m+\underline{u}_{\max}}=\lim
\inf_{m\rightarrow E^{-}}\Psi^{-1}(m+\underline{u}_{\max})\geq\lim
\inf_{m\rightarrow E^{-}}\Psi^{-1}(m)$ $=\infty$, we obtain $\lim
_{m\rightarrow E^{-}}c_{m+\underline{u}_{\max}}=\infty$. Therefore, using
Lemmas 3.4 and 3.9, we have $\lim_{m\rightarrow E^{-}}\underline{u}_{\max}$
$=\lim_{c\rightarrow\infty}c\Psi^{\prime}(c)$, provided one of them exists.
The rest of this lemma is deduced from Lemma 3.10.\hfill$\square$

L{\small EMMA} 3.12 [10, Lemma 1.1.2]. \textit{If }$\lim_{c\rightarrow
(-\xi)^{+}}\Psi^{\prime}(c)=\infty,$\textit{ \ }$\lim_{c\rightarrow(-\xi)^{+}%
}(-1)^{n-1}$ $\Psi^{(n)}(c)$ $=\infty$\textit{ }$(n=1,2,3,...)$.

P{\small ROOF}. From (3.3), we have $(-1)^{n-1}\Psi^{(n)}(c)$ $=n!\int
_{0}^{\infty}(t+c+\xi)^{-n-1}\rho(dt)\geq0$ $(c>-\xi$, $n=1,2,3,...)$. Since
$\Psi^{(3)}(c)\geq0,$ we observe that $\Psi^{\prime}(a)-\Psi^{\prime}(c)$
$=\int_{c}^{a}\Psi^{(2)}(x)dx\geq\Psi^{(2)}(c)(a-c)$ for each $-\xi<c<a$. This
implies that $\Psi^{(2)}(c)$ $\leq\left(  \Psi^{\prime}(a)-\Psi^{\prime
}(c)\right)  /(a-c)$. Thus, using $\lim_{c\rightarrow(-\xi)^{+}}\Psi^{\prime
}(c)$ $=\infty$, we have $\lim_{c\rightarrow(-\xi)^{+}}\Psi^{(2)}(c)$
$=-\infty$. For each $n\in\{3,4,5,...\},$ we find that $(-1)^{n}(\Psi
^{(n-1)}(a)-\Psi^{(n-1)}(c))$ $=\int_{c}^{a}(-1)^{n}\Psi^{(n)}(x)dx$
$\geq(-1)^{n}\Psi^{(n)}(c)(a-c)$, which implies that $(-1)^{n-1}\Psi^{(n)}(c)$
$\geq(-1)^{n-2}\Psi^{(n-1)}(c)/(a-c)$ $+(-1)^{n-1}\Psi^{(n-1)}(a)/(a-c).$
Therefore, by induction on $n$, we arrive at the conclusion.\hfill$\square$

L{\small EMMA} 3.13. $\lim_{c\rightarrow(-\xi)^{+}}(c+\xi)^{n+1}\Psi
^{(n)}(c)=0$ $(n=1,2,3,...)$.

P{\small ROOF}. From (3.3), we have $(c+\xi)^{n+1}\Psi^{(n)}(c)$
$=(-1)^{n-1}n!\int_{0}^{\infty}(c+\xi)^{n+1}/(t$ $+c+\xi)^{n+1}\rho(dt)$.
Thus, applying Lebesgue's (monotone convergence) theorem, we arrive at the
conclusion.\hfill$\square$

T{\small HEOREM} 3.14. \textit{If }$H_{\xi}=\infty,$\textit{ then }%
$\lim_{m\rightarrow(\xi+1/H_{\xi})^{+}}\underline{u}_{\max}$\textit{ }%
$=0.$\textit{ If }$H_{\xi}<\infty,$\textit{ then }$\lim_{m\rightarrow
(\xi+1/H_{\xi})^{+}}\underline{u}_{\max}$\textit{ }$>0$\textit{, which can
assume any positive value, exists.}

P{\small ROOF}. From (3.3), we have $\Psi^{\prime}(c)>0$ and $\Psi
^{\prime\prime}(c)<0$ for $c>-\xi$. As shown in the proof of Lemma 2.5,
$c_{m+\underline{u}_{\max}}^{\prime}=1/\Psi^{\prime}(m+\underline{u}_{\max
})>0$. Therefore, $\lambda$ $:=\lim_{m\rightarrow(\xi+1/H_{\xi})^{+}%
}c_{m+\underline{u}_{\max}}$ exists such that $-\xi\leq\lambda<\infty$. Since
$\underline{u}_{\max}=\Psi(c_{m+\underline{u}_{\max}})$ $-m$, $\lim
_{m\rightarrow(\xi+1/H_{\xi})^{+}}\underline{u}_{\max}$ $=\Psi(\lambda)$
$-\xi$ $-1/H_{\xi}$. From $\Psi^{\prime\prime}(c)<0$, $\Psi^{\prime}%
(\lambda^{+})>0$ exists (including $+\infty$). From the proof of Lemma 3.11,
we observe that $\underline{u}_{\max}$ $=\left(  c_{m+\underline{u}_{\max}%
}+\Psi(c_{m+\underline{u}_{\max}})\right)  $ $/\left(  1+1/\Psi^{\prime
}(c_{m+\underline{u}_{\max}})\right)  $, which induces $\Psi(\lambda
)-\xi-1/H_{\xi}$ $=\left(  \lambda+\Psi(\lambda)\right)  $ $/\left(
1+1/\Psi^{\prime}(\lambda^{+})\right)  .$

If $H_{\xi}=\infty$, then $\Psi(-\xi)=\xi.$ For each $c>-\xi$, put $m:=\left(
\Psi(c)-c\Psi^{\prime}(c)\right)  $ $/(1+\Psi^{\prime}(c))$ and $v:=\Psi
(c)-m$. Then%
\[
\left\{
\begin{array}
[c]{c}%
m+v=\Psi(c),\\
v=(m+c)\Psi^{\prime}(c).
\end{array}
\right.
\]
Moreover, from $\partial m/\partial c=-(c+\Psi(c))\Psi^{\prime\prime
}(c)/(1+\Psi^{\prime}(c))^{2}>0,$ we obtain $\xi<m$ $<E$. Thus, in accordance
with (2.2), we can consider $v=\underline{u}_{\max}$ and $c=c_{m+\underline
{u}_{\max}}$. Therefore, $-\xi\leq\lambda\leq c_{m+\underline{u}_{\max}}=c$
for each $c>-\xi$, which implies that $\lambda=-\xi$. Thus, $\lim
_{m\rightarrow(\xi+1/H_{\xi})^{+}}\underline{u}_{\max}$ $=\Psi(-\xi)-\xi=0$.

If $H_{\xi}<\infty$, then from $\left(  \lambda+\Psi(\lambda)\right)  /\left(
1+1/\Psi^{\prime}(\lambda^{+})\right)  $ $\geq1/(H_{\xi}\left(  1+1/\Psi
^{\prime}(\lambda^{+})\right)  )$ $>0$, $\lim_{m\rightarrow(\xi+1/H_{\xi}%
)^{+}}\underline{u}_{\max}>0$. In Example 2.7, we observe that $\xi+1/H_{\xi
}=r/3$ and $\underline{u}_{\max}$ $=(3r-m)(m+r)/(m+9r)$. Thus, we obtain
$\lim_{m\rightarrow(r/3)^{+}}\underline{u}_{\max}=8r/21$ $(r>0)$, which
implies the conclusion.\hfill$\square$

The following Lemma is similar to [11, Lemma 2.10].

L{\small EMMA} 3.15. $\left\vert \Psi^{(n+1)}(c)/\Psi^{(n)}(c)\right\vert
<(n+1)/(c+\xi)$ \textit{and }$\Psi^{\prime}(c)/\left(  \Psi(c)-\xi-1/H_{\xi
}\right)  $ $<1/(c+\xi)$ $(c>-\xi$, $n=1,2,3,...)$.

P{\small ROOF}. From (3.3), we have
\begin{align*}
\left\vert \Psi^{(n+1)}(c)\right\vert  &  =(n+1)!\int_{0}^{\infty}\frac
{1}{(t+c+\xi)^{n+2}}\rho(dt)\\
&  =n!\int_{0}^{\infty}\frac{n+1}{t+c+\xi}\times\frac{1}{(t+c+\xi)^{n+1}}%
\rho(dt)\\
&  <\frac{n+1}{c+\xi}\times n!\int_{0}^{\infty}\frac{1}{(t+c+\xi)^{n+1}}%
\rho(dt)\\
&  =\frac{n+1}{c+\xi}\left\vert \Psi^{(n)}(c)\right\vert \text{ }%
(n=1,2,3,...).
\end{align*}
Moreover, from (3.1), we observe that%
\[
(c+\xi)\Psi^{\prime}(c)=\int_{0}^{\infty}\frac{c+\xi}{t(t+c+\xi)}\rho
(dt)-\int_{0}^{\infty}\frac{(c+\xi)^{2}}{t(t+c+\xi)^{2}}\rho(dt)<\Psi
(c)-\xi-1/H_{\xi}.
\]

\hfill$\square$

C{\small ORORALLY} 3.16. $\lim\inf_{m\rightarrow E^{-}}\underline{u}_{\max
}^{\prime}\geq-1/2.$

P{\small ROOF}. From Lemma 2.5, we have $\underline{u}_{\max}%
=(m+c_{m+\underline{u}_{\max}})\Psi^{\prime}(c_{m+\underline{u}_{\max}})$,
$\ \Psi^{\prime}(c_{m+\underline{u}_{\max}})$ $c_{m+\underline{u}_{\max}%
}^{\prime}$ $=1$ and $m=\left(  \Psi(c_{m+\underline{u}_{\max}}%
)-c_{m+\underline{u}_{\max}}\Psi^{\prime}(c_{m+\underline{u}_{\max}})\right)
/(1+\Psi^{\prime}(c_{m+\underline{u}_{\max}}))$. It follows that
\begin{align*}
\underline{u}_{\max}^{\prime}  &  =-1-\frac{1+\Psi^{\prime}(c_{m+\underline
{u}_{\max}})}{(m+c_{m+\underline{u}_{\max}})c_{m+\underline{u}_{\max}}%
^{\prime}\Psi^{\prime\prime}(c_{m+\underline{u}_{\max}})}\\
&  =-1-\frac{\Psi^{\prime}(c_{m+\underline{u}_{\max}})}{(c_{m+\underline
{u}_{\max}}+\xi)\Psi^{\prime\prime}(c_{m+\underline{u}_{\max}})}\times
\frac{(c_{m+\underline{u}_{\max}}+\xi)(1+\Psi^{\prime}(c_{m+\underline
{u}_{\max}}))^{2}}{c_{m+\underline{u}_{\max}}(1+\frac{\Psi(c_{m+\underline
{u}_{\max}})}{c_{m+\underline{u}_{\max}}})}.
\end{align*}
On the other hand, from to Lemmas 2.1, 3.4, and 3.9, we have%
\[
\lim_{m\rightarrow E^{-}}\frac{c_{m+\underline{u}_{\max}}+\xi}{c_{m+\underline
{u}_{\max}}}\cdot\frac{(1+\Psi^{\prime}(c_{m+\underline{u}_{\max}}))^{2}%
}{1+\frac{\Psi(c_{m+\underline{u}_{\max}})}{c_{m+\underline{u}_{\max}}}}=1.
\]
From Lemma 3.15, we observe that $-\Psi^{\prime\prime}(c)/\Psi^{\prime
}(c)<2/(c+\xi),$ which implies that $-1-\Psi^{\prime}(c)/\left(  (c+\xi
)\Psi^{\prime\prime}(c)\right)  >-1/2$. Therefore, we conclude that%
\[
\lim\inf_{m\rightarrow E^{-}}\underline{u}_{\max}^{\prime}=\lim\inf
_{c\rightarrow\infty}\left(  -1-\frac{\Psi^{\prime}(c)}{\left(  c+\xi\right)
\Psi^{\prime\prime}(c)}\right)  \geq-\frac{1}{2}.
\]

\hfill$\square$

L{\small EMMA} 3.17. \textit{Assume} $E<\infty$. \textit{Then},%
\begin{equation}
\lim_{m\rightarrow E^{-}}\frac{\underline{u}_{\max}}{\underline{E}}=\frac
{1}{\frac{1}{\lim_{c\rightarrow\infty}\frac{c\Psi^{\prime}(c)}{E-\Psi(c)}}%
+1}=\frac{1}{\frac{1}{\lim_{c\rightarrow\infty}\frac{\int_{0}^{\infty}\left(
c/(c+t)\right)  ^{2}\rho(dt)}{\int_{0}^{\infty}c/(c+t)\rho(dt)}}+1} \tag{3.5}%
\end{equation}
\textit{if one of three limits exists.}

P{\small ROOF}. From (2.2) and (3.4), we have%
\[
\frac{\underline{u}_{\max}}{\underline{E}}=\frac{1}{\frac{E-\Psi
(c_{m+\underline{u}_{\max}})}{c_{m+\underline{u}_{\max}}\Psi^{\prime
}(c_{m+\underline{u}_{\max}})\left(  1+\Psi(c_{m+\underline{u}_{\max}%
})/c_{m+\underline{u}_{\max}}\right)  }+\frac{c_{m+\underline{u}_{\max}}%
+E}{c_{m+\underline{u}_{\max}}+\Psi(c_{m+\underline{u}_{\max}})}}.
\]
Thus, from Lemmas 1.1, 2.1, 3.4, and 3.10, we obtain $\lim_{m\rightarrow
E^{-}}\underline{u}_{\max}/\underline{E}$ $=1/(1/\lim_{c\rightarrow\infty
}(c\Psi^{\prime}(c)$ $/(E-\Psi(c)))+1)$. Using (3.1) and (3.3) we observe that%
\begin{equation}
\frac{c\Psi^{\prime}(c)}{E-\Psi(c)}=\frac{\int_{0}^{\infty}\left(
(c+\xi)/(t+c+\xi)\right)  ^{2}\rho(dt)}{(1+\frac{\xi}{c})\int_{0}^{\infty
}(c+\xi)/(t+c+\xi)\rho(dt)}, \tag{3.6}%
\end{equation}
which yields the desired equation.\hfill$\square$

L{\small EMMA} 3.18.
\begin{equation}
-\lim_{m\rightarrow E^{-}}\underline{u}_{\max}^{\prime}=1+\lim_{c\rightarrow
\infty}\frac{\Psi^{\prime}(c)}{\left(  c+\xi\right)  \Psi^{\prime\prime}%
(c)}=1-\lim_{c\rightarrow\infty}\frac{\int_{0}^{\infty}(c/(c+t))^{2}\rho
(dt)}{2\int_{0}^{\infty}(c/(c+t))^{3}\rho(dt)}, \tag{3.7}%
\end{equation}
\textit{if one of three limits exists. In this case, if }$E<\infty$\textit{,
its value is equal to }$\lim_{m\rightarrow E^{-}}\underline{u}_{\max}$
$/\underline{E}$.

P{\small ROOF}. From (3.3), we have $\Psi^{\prime}(c)/(\left(  c+\xi\right)
\Psi^{\prime\prime}(c))$ $=-\int_{0}^{\infty}(\left(  c+\xi\right)
/(t+c+\xi))^{2}\rho(dt)$ $/\left(  2\int_{0}^{\infty}(\left(  c+\xi\right)
/(t+c+\xi))^{3}\rho(dt)\right)  $. From the proof of Corollary 3.16, we have
$\lim_{m\rightarrow E^{-}}\underline{u}_{\max}^{\prime}$ $=-1-\lim
_{c\rightarrow\infty}\Psi^{\prime}(c)/(\left(  c+\xi\right)  \Psi
^{\prime\prime}(c))$ if one of them exists. In this case, if $E<\infty,$ then
from $\lim_{m\rightarrow E^{-}}\underline{u}_{\max}=0$ (Lemma 3.11) and by
using the mean value theorem, we obtain $\lim_{m\rightarrow E^{-}}%
\underline{u}_{\max}/\underline{E}$ $=-\lim_{m\rightarrow E^{-}}\int_{m}%
^{E}\underline{u}_{\max}^{\prime}(t)dt/\underline{E}$ $=-\lim_{m\rightarrow
E^{-}}\underline{u}_{\max}^{\prime}$, which implies the conclusion.\hfill
$\square$

T{\small HEOREM} 3.19.\textbf{ }\textit{If }$\int_{0}^{\infty}\rho(dt)<\infty
$, $\lim_{m\rightarrow E^{-}}\underline{u}_{\max}/\underline{E}=1/2$.

P{\small ROOF}. Using $\int_{0}^{\infty}\rho(dt)<\infty$ and $\int_{0}%
^{\infty}1/\left(  t(t+1)\right)  \rho(dt)<\infty$, we have $E<\infty$ (Lemma
3.6). By applying Lebesgue's theorem, we observe that $\lim_{c\rightarrow
\infty}(\int_{0}^{\infty}c/(c+t)\rho(dt)$ $/\int_{0}^{\infty}\left(
c/(c+t)\right)  ^{2}\rho(dt))$ $=1.$ Thus, from Lemma 3.17, we obtain
$\lim_{m\rightarrow E^{-}}\underline{u}_{\max}$ $/\underline{E}=1/2$%
.\hfill$\square$

C{\small ORORALLY} 3.20.\textit{ }$\lim_{m\rightarrow E^{-}}\underline
{u}_{\max}/\underline{E}$\textit{ }$=1/2$\textit{ for any finite game
}$\{(a_{j}$\textit{, }$p_{j})\}$\textit{ such that }$\sum_{j=1}^{n}a_{j}%
p_{j}>0,$\textit{ }$\sum_{j=1}^{n}p_{j}=1$, $0\leq p_{j}<1$\textit{, and
}$1\leq j\leq n$.

P{\small ROOF}.\textbf{ }The complete Bernstein function
\[
\Psi(c)-\xi-\frac{1}{H_{\xi}}=\frac{1}{\sum_{j=1}^{n}\frac{p_{j}}{a_{j}+c}%
}-c-\xi-1/H_{\xi}\text{ \ \ \ }(c>-\xi)
\]
is a rational function with respect to $c$, which is analytic on $%
\mathbb{C}
\backslash(-\infty,$ $-\xi]$ and preserves the upper and lower half-planes
(Theorems 3.1 and 3.3). Therefore, using [8, Theorem 2.2 (vi)], we obtain the
representation $\Psi(c)-\xi-1/H_{\xi}$ $=-\sum_{j=1}^{m}e_{j}/(c$ $+d_{j})$
$+$ $k$, where $e_{j}>0$, $d_{j}>\xi$ and $k$ is a constant. Since $\Psi
(-\xi)=\xi+1/H_{\xi}$, $k$ $=\sum_{j=1}^{m}e_{j}/(d_{j}$ $-\xi)$ holds.
Defining $\sigma(dt)$\ as the sum of Dirac measures $\sum_{j=1}^{m}%
e_{j}/(d_{j}$ $-\xi)\delta_{d_{j}-\xi}$, we observe that
\[
\int_{0}^{\infty}\frac{c+\xi}{t+c+\xi}\sigma(dt)=-\sum_{j=1}^{m}\frac{e_{j}%
}{c+d_{j}}+\sum_{j=1}^{m}\frac{e_{j}}{d_{j}-\xi}=\Psi(c)-\xi-\frac{1}{H_{\xi}%
}.
\]
Using Theorem 3.1, we obtain $\int_{0}^{\infty}\rho(dt)=\int_{0}^{\infty
}t\sigma(dt)=\sum_{j=1}^{m}e_{j}<\infty$, which, in accordance with Theorem
3.19, implies the conclusion.\hfill$\square$

\bigskip

\begin{center}
\textbf{4. Abelian theorems\smallskip}
\end{center}

In the following paragraphs, we assume that a nonzero measure $\rho(dt)$
originates from (3.1). For a function $f(x)>0$, $\omega_{f}:=\lim
\sup_{x\rightarrow\infty}$ $\log f(x)/\log x$ is termed \textit{the upper
order }(see [1, Section 2.2.2]). We will show that $\lim_{m\rightarrow E^{-}%
}\underline{u}_{\max}/\underline{E}$ can be calculated by the upper order of
the function $\int_{0}^{x}\rho(dt)$.

A measurable function $f(x)>0$ is said to be \textit{regularly varying }of
\textit{index} $r$, written as $f\in R_{r}$, if $\lim_{x\rightarrow\infty}$
$f(\lambda x)/f(x)=\lambda^{r}$ for each $\lambda>0$ (see [1, Section 1.4.2]).
It is easy to verify that $\omega_{f}=r$ if $f\in R_{r}$. The notation $l(x)$
is used only for a \textit{slowly varying} function such that $l(x)\in R_{0}$.
We write $f(x)\sim cg(x)$ when $\lim_{x\rightarrow\infty}f(x)/g(x)=c$. If
$c=0$, the relation $f(x)\sim cg(x)$ suggests that $f(x)=o(g(x))$ (see [1, Preface]).

L{\small EMMA} 4.1. \textit{If }$\int_{0}^{x}\rho(dt)\in R_{r}$\textit{,
}$0\leq r\leq2$\textit{. In this case, if }$E<\infty,$ $0\leq r\leq1$.

P{\small ROOF}. We can write $\int_{0}^{x}\rho(dt)=x^{r}l(x)$ with $l(x)\in
R_{0}$. Assuming $r<0$, then $\lim_{x\rightarrow\infty}\int_{0}^{x}%
\rho(dt)=\lim_{x\rightarrow\infty}x^{r}l(x)=0$ (see [1, Proposition 1.3.6]),
which contradicts the fact that $\int_{0}^{\infty}\rho(dt)>0$. Assuming $r>2$,
we have $\lim_{x\rightarrow\infty}\int_{0}^{x}\rho(dt)$ $/(x(x+1))$
$=\lim_{x\rightarrow\infty}x^{r-2}l(x)$ $\cdot\lim_{x\rightarrow\infty
}1/(1+1/x)$ $=\infty$. On the other hand, for each $x>0$, we have $\int
_{0}^{x}\rho(dt)/(x(x+1))\leq\int_{0}^{\infty}(t(t+1))^{-1}\rho(dt)<\infty$
(Lemma 3.5), which is a contradiction.

If $E<\infty$, then $\int_{0}^{\infty}t^{-1}\rho(dt)<\infty$ (Lemma 3.6).
Thus, for each $x>0$, we have $\int_{0}^{x}\rho(dt)/x\leq\int_{0}^{\infty
}t^{-1}\rho(dt)<\infty$, which implies that $r\leq1$.\hfill$\square$

L{\small EMMA} 4.2.\textbf{ }\textit{Suppose }$0\leq r<n$\textit{. Then,
}$\int_{0}^{x}\rho(dt)\sim x^{r}l(x)$\textit{ }$(x\rightarrow\infty)$\textit{
if and only if} $\int_{0}^{\infty}(x/(x+t))^{n}\rho(dt)\sim\Gamma
(n-r)\Gamma(r+1)x^{r}l(x)/\Gamma(n)$ $(x\rightarrow\infty)$.

P{\small ROOF}.\textbf{ }The nondecreasing function\textbf{ }$U(x):=\int
_{0}^{x}\rho(dt)$ satisfies $U(0^{-})=0.$ Since $0<n-r\leq n$, using [1,
Theorem 1.7.4], we obtain that $U(x)\sim x^{r}l(x)$ $(x\rightarrow\infty)$ is
equivalent to%

\[
\int_{0}^{\infty}\frac{d(U(t))}{(x+t)^{n}}\sim\frac{\Gamma(n-r)\Gamma
(r+1)}{\Gamma(n)}x^{r-n}l(x)\text{ \ \ \ }(x\rightarrow\infty),
\]
which implies the conclusion.\hfill$\square$

L{\small EMMA} 4.3. \textit{If }$\int_{0}^{x}\rho(dt)\in R_{r}$\textit{ and
}$r\neq2$\textit{, then }$\lim_{c\rightarrow\infty}\Psi^{\prime}(c)/\left(
c\Psi^{\prime\prime}(c)\right)  $ $=1/\left(  r-2\right)  $.

P{\small ROOF}. From Lemma 4.1, we have $0\leq r<2$. From Lemma 4.2, we
observe that%
\[
\lim_{x\rightarrow\infty}\frac{\int_{0}^{\infty}(x/(x+t))^{2}\rho(dt)}%
{\int_{0}^{\infty}(x/(x+t))^{3}\rho(dt)}=\frac{\Gamma(3)\Gamma(2-r)}%
{\Gamma(2)\Gamma(3-r)}=\frac{2}{2-r}.
\]
Therefore, the relation (3.7) implies the conclusion. \hfill$\square$

T{\small HEOREM} 4.4. \textit{If }$E<\infty$\textit{ and }$\int_{0}^{x}%
\rho(dt)\in R_{r}$\textit{, then }$\lim_{m\rightarrow E^{-}}\underline
{u}_{\max}/\underline{E}$\textit{ }$=\left(  1-r\right)  /\left(  2-r\right)
$\textit{ with }$0\leq r\leq1$\textit{.}

P{\small ROOF}. From Lemmas 3.18 and 4.3, we obtain $\lim_{m\rightarrow E^{-}%
}\underline{u}_{\max}/\underline{E}$ $=1+\lim_{c\rightarrow\infty}\Psi
^{\prime}(c)$ $/\left(  (c+\xi)\Psi^{\prime\prime}(c)\right)  =(1-r)/\left(
2-r\right)  $.\hfill$\square$

Whenever we use the notation $q(t)$, it is understood that $\rho(dt)=q(t)dt$
with $q(t)\geq0$.

L{\small EMMA} 4.5. \textit{If }$q(t)\in R_{\alpha}$\textit{, then }%
$\alpha\leq1$\textit{. In addition, if }$E<\infty,$\textit{ then }$\alpha
\leq0$.

P{\small ROOF}. We can write $q(t)=t^{\alpha}l(t)$ with $l(t)\in R_{0}$. From
[1, Corollary 1.4.2], $X>0$ exists such that $l(x)$ is locally bounded in
$[X$, $\infty)$. Assuming $\alpha>1$, then using [1, Propositions 1.3.6 and
1.5.8], we obtain $\lim_{x\rightarrow\infty}\int_{X}^{x}\rho(dt)/(x(x+1))$
$=\lim_{x\rightarrow\infty}x^{\alpha-1}l(x)/(\alpha+1)$ $\cdot\lim
_{x\rightarrow\infty}1/(1+1/x)=\infty$. On the other hand, for each $x>X$, we
observe that $\int_{X}^{x}\rho(dt)/(x(x+1))$ $\leq\int_{X}^{x}(t(t+1))^{-1}%
\rho(dt)$ $\leq\int_{0}^{\infty}(t(t$ $+1))^{-1}\rho(dt)<\infty$ (Lemma 3.5),
which is a contradiction.

When $E<\infty$, we have $\int_{X}^{x}\rho(dt)/x\leq\int_{0}^{\infty}%
t^{-1}\rho(dt)<\infty$ (Lemma 3.6). Thus, arguments similar to the one above
yield the conclusion.\hfill$\square$

L{\small EMMA} 4.6. \textit{If }$q(t)\in R_{\alpha}$\textit{, then }$\int
_{0}^{x}q(t)dt\in R_{\alpha+1}$\textit{ }$(\alpha>-1)$\textit{ or }$\int
_{0}^{x}q(t)dt\in R_{0}$ $(\alpha\leq-1)$.

P{\small ROOF}. As the proof of Lemma 4.5, if $\alpha>-1$, using [1,
Proposition 1.5.8], we obtain $\int_{X}^{x}t^{\alpha}l(t)dt\sim x^{\alpha
+1}l(x)/(\alpha+1)\in R_{\alpha+1}$. If $\alpha=-1$, then from [1, Proposition
1.5.9a], we have $\int_{X}^{x}t^{-1}l(t)dt\in R_{0}$. If $\alpha<-1$, then the
nondecreasing function $\int_{0}^{x}q(t)dt=\int_{0}^{x}t^{\alpha}l(t)dt$ is
bounded as will be shown below, which suggests that $\int_{0}^{x}q(t)dt\in
R_{0}$. Put $\varepsilon:=-(\alpha+1)/2>0$. Then, from $\lim_{t\rightarrow
\infty}t^{-\varepsilon}l(t)$ $=0$, $Y>0$ exists such that $0\leq
t^{-\varepsilon}l(t)\leq1$ for each $t\geq Y$. Therefore, for each $x\geq Y$,
we find that $\int_{0}^{x}q(t)dt\leq\int_{0}^{Y}q(t)dt$ $+\int_{Y}%
^{x}t^{-\varepsilon}l(t)t^{-1-\varepsilon}dt$ $\leq\int_{0}^{Y}%
q(t)dt+1/(\varepsilon Y^{\varepsilon})$ $<\infty$.

\hfill$\square$

C{\small ORORALLY} 4.7.\textit{ If }$E<\infty$\textit{ and }$q(t)\in
R_{\alpha}$\textit{, then }$\alpha\leq0$\textit{ and}$\ $%
\begin{equation}
\lim_{m\rightarrow E^{-}}\frac{\underline{u}_{\max}}{\underline{E}}=\left\{
\begin{array}
[c]{cl}%
\frac{1}{2}, & \text{\textit{if} \ \ \ \ }\alpha\leq-1,\\
\frac{\alpha}{\alpha-1}, & \text{\textit{if} }-1<\alpha\leq0.
\end{array}
\right.  \tag{4.1}%
\end{equation}

P{\small ROOF}. It is the direct consequence of Theorem 4.4 and Lemmas 4.5 and 4.6.

\hfill$\square$

\newpage

\begin{center}
\textbf{5. Tauberian theorems}
\end{center}

Given a measurable function $f:(0,\infty)\rightarrow%
\mathbb{R}
,$ let $\check{f}(z):=\int_{0}^{\infty}t^{-z-1}f(t)dt$ be its Mellin transform
for $z\in%
\mathbb{C}
$ such that the integral converges absolutely (see [1, 2, 3]). For example,
putting $k(x):=2x^{2}$ $(0<x<1)$ or $0$ $(x\geq1)$, we obtain $\check
{k}(z)=2/(2-z)$ $(-\infty<\operatorname{Re}z<2)$. In addition, putting
$h(x):=x$ $(0<x<1)$ or $0$ $(x\geq1)$, we obtain $\check{h}(z)=1/(1-z)$
$(-\infty<\operatorname{Re}z<1)$. Given measurable functions $f,$
$g:(0,\infty)\rightarrow%
\mathbb{R}
,$ let $(f\ast g)(x)$ $:=\int_{0}^{\infty}t^{-1}f(x/t)g(t)dt$ be the Mellin
convolution of these functions for $x>0$ such that the integral converges absolutely.

T{\small HEOREM} 5.1.\textbf{ }\textit{If }$E<\infty$\textit{ and }%
$\lim_{m\rightarrow E^{-}}\underline{u}_{\max}^{\prime}$\textit{ exists, then
}$\lim_{m\rightarrow E^{-}}\underline{u}_{\max}/\underline{E}$\textit{
}$=(1-r)/(2-r)$\textit{ }$(0\leq r\leq1)$\textit{ and }$\int_{0}^{x}%
\rho(dt)\in R_{r}$\textit{, where }$r$\textit{ is the upper order of }%
$\int_{0}^{\infty}\left(  x/(x+t)\right)  ^{n}\rho(dt)$\textit{ }%
$(n>1)$\textit{.}

P{\small ROOF}. Putting $K(x):=\int_{0}^{\infty}\left(  x/(x+t)\right)
^{3}\rho(dt)$, we observe that $\check{k}(2^{-})=\infty$ and%
\[
\frac{\left(  k\ast K\right)  (x)}{K(x)}=\frac{\int_{0}^{\infty}\left(
x/(x+t)\right)  ^{2}\rho(dt)}{\int_{0}^{\infty}\left(  x/(x+t)\right)
^{3}\rho(dt)}\geq1.
\]
From Lemma 3.18, $\lim_{x\rightarrow\infty}\left(  k\ast K\right)
(x)/K(x)=c\geq1$ exists. As $K(x)$ is an increasing function, from [1, Theorem
5.2.3 and Section 2.1.2], we obtain $c=\check{k}(\omega_{K})=2/(2-\omega_{K}%
)$, $\omega_{K}<2$ and $K(x)$ is regularly varying. Thus, $\lim_{m\rightarrow
E^{-}}\underline{u}_{\max}/\underline{E}$ $=-\lim_{m\rightarrow E^{-}%
}\underline{u}_{\max}^{\prime}$\ $=1-1/(2-\omega_{K})$ $=(1-\omega
_{K})/(2-\omega_{K})$. From Lemma 4.2, it follows that $\int_{0}^{x}%
\rho(dt)\in R_{\omega_{K}}$. Moreover, from Lemma 4.1, we have $0\leq
\omega_{K}\leq1$. This implies that the upper order of $\int_{0}^{\infty
}\left(  x/(x+t)\right)  ^{n}\rho(dt)$ is always $\omega_{K}$ for each
$n>1$.\hfill$\square$

T{\small HEOREM} 5.2.\textbf{ }\textit{If }$E<\infty$\textit{ and }%
$\lim_{m\rightarrow E^{-}}\underline{u}_{\max}/\underline{E}\neq0$\textit{
exists, then }$\lim_{m\rightarrow E^{-}}\underline{u}_{\max}$ $/\underline
{E}=-\lim_{m\rightarrow E^{-}}\underline{u}_{\max}^{\prime}$\textit{
}$=(1-r)/(2-r)$\textit{ }$(0\leq r<1)$\textit{ and }$\int_{0}^{x}\rho(dt)\in
R_{r}$\textit{, where }$r$\textit{ is the upper order of }$\int_{0}^{\infty
}\left(  x/(x+t)\right)  ^{n}\rho(dt)$\textit{ }$(n>1)$\textit{.}

P{\small ROOF}. This proof is formally the same as that in Theorem 5.1.
Putting $S(x):=\int_{0}^{\infty}\left(  x/(x+t)\right)  ^{2}\rho(dt)$, we
observe that $\check{h}(1^{-})=\infty$ and%
\[
\frac{\left(  h\ast S\right)  (x)}{S(x)}=\frac{\int_{0}^{\infty}%
x/(x+t)\rho(dt)}{\int_{0}^{\infty}\left(  x/(x+t)\right)  ^{2}\rho(dt)}\geq1.
\]
From Lemma 3.17, $\lim_{x\rightarrow\infty}\left(  h\ast S\right)
(x)/S(x)=c\geq1$ exists. As $S(x)$ is an increasing function, from [1, Theorem
5.2.3 and Section 2.1.2], we obtain $c=\check{h}(\omega_{S})$ $=1/(1$
$-\omega_{S})$, $\omega_{S}<1,$\ and that $S(x)$ is regularly varying. Thus,
$\lim_{m\rightarrow E^{-}}\underline{u}_{\max}/\underline{E}$ $=1/(c+1)$
$=(1-\omega_{S})/(2-\omega_{S})$. From Lemma 4.2, it follows that $\int
_{0}^{x}\rho(dt)\in R_{\omega_{S}}$. Moreover, from Lemma 4.1, we have
$0\leq\omega_{S}<1$. This implies that the upper order of $\int_{0}^{\infty
}\left(  x/(x+t)\right)  ^{n}\rho(dt)$ is always $\omega_{K}$ for each $n>1$.
In this case, based on Lemmas 3.18 and 4.3, $\lim_{m\rightarrow E^{-}%
}\underline{u}_{\max}^{\prime}$\ exists.\hfill$\square$

C{\small ORORALLY} 5.3.\textbf{ }\textit{If }$E<\infty$\textit{, the following
are equivalent.}

(a)\textit{ }$\lim_{m\rightarrow E^{-}}\underline{u}_{\max}/\underline{E}%
\neq0$\textit{ exists.}

(b) $\int_{0}^{x}\rho(dt)\in R_{r}$ $(r\neq1).$

(c)\textit{ }$\lim_{m\rightarrow E^{-}}\underline{u}_{\max}^{\prime}\neq
0$\textit{ exists.}

P{\small ROOF}. Using Theorems 4.4, 5.1, and 5.2, and Lemmas 3.18 and 4.3, we
arrive at the conclusion.\hfill$\square$

\bigskip

It is noteworthy that $\int_{0}^{x}\rho(dt)\in R_{1}$ includes $\lim
_{m\rightarrow E^{-}}\underline{u}_{\max}/\underline{E}=0$. However, the
converse is not necessarily true because a nonregularly varying function
$\int_{e}^{x}\rho(dt)$ $=(2+\sin\left(  \log x\right)  )x/(1+\log
x)^{3/2}-(2+\sin1)e/(2\sqrt{2})$ $(x\geq e)$ provides an example with
$\lim_{m\rightarrow E^{-}}\underline{u}_{\max}/\underline{E}=0$. The details
are left to the reader. In this direction, we observe that $\lim_{m\rightarrow
E^{-}}\underline{u}_{\max}/\underline{E}=0$ if and only if $\int_{0}^{\infty
}(x+t)^{-1}\rho(dt)$ is normalized slowly varying. Because, since
$\lim_{c\rightarrow\infty}c\Psi^{\prime}(c)/(E$ $-\Psi(c))=0$ (Lemma 3.17),
$E-\Psi(c)$ $=\int_{0}^{\infty}(t+c+\xi)^{-1}\rho(dt)$ (Lemma 3.5) is
normalized slowly varying (see [1, (1.3.4)]).

L{\small EMMA} 5.4 [1, Theorem 1.7.2]\textbf{.} \textit{If }$\int_{0}%
^{x}f(t)dt\sim cx^{r}l(x)$\textit{ }$(x\rightarrow\infty)$\textit{, where
}$f(x)$\textit{ is nondecreasing or nonincreasing in an interval }%
$(T$\textit{, }$\infty)$\textit{ }$(T>0)$\textit{, then }$f(x)$ $\sim
crx^{r-1}l(x)$\textit{ }$(x\rightarrow\infty)$\textit{.}

C{\small ORORALLY} 5.5.\textbf{ }\textit{Assuming }$E<\infty,$\textit{ }%
$\rho(dt)=q(t)dt$\textit{ and }$q(t)$\textit{ is nonincreasing in an interval
}$(T$\textit{, }$\infty)$\textit{ }$(T>0)$\textit{. When }$\lim_{m\rightarrow
E^{-}}\underline{u}_{\max}/\underline{E}\neq0$\textit{ exists, the following
properties hold:}

(1) \textit{If }$\omega_{q}\leq-1$\textit{, then }$\lim_{m\rightarrow E^{-}%
}\underline{u}_{\max}/\underline{E}=1/2$\textit{ and }$\int_{0}^{x}%
q(t)dt$\textit{ is slowly varying.}

(2) \textit{If} $\omega_{q}>-1$, \textit{then} $\omega_{q}<0$, $\lim
_{m\rightarrow E^{-}}\underline{u}_{\max}/\underline{E}=-\omega_{q}%
/(1-\omega_{q})$, \textit{and} $q(t)$ \textit{is regularly varying}.

P{\small ROOF}. Put $S(x):=\int_{0}^{\infty}\left(  x/(x+t)\right)
^{2}q(t)dt.$ Then, by applying Theorem 5.2, we obtain $\lim_{m\rightarrow
E^{-}}\underline{u}_{\max}/\underline{E}$ $=(1-\omega_{S})/(2-\omega_{S}),$
\textit{ }$(0\leq\omega_{S}<1)$ and $\int_{0}^{x}q(t)dt$ $\in R_{\omega_{S}}$.
Thus, by Lemma 5.4 we find that $\int_{0}^{x}q(t)dt\sim x^{\omega_{S}}%
l(x)$\textit{ }$(x\rightarrow\infty)$ and $q(t)$ $\sim\omega_{S}t^{\omega
_{S}-1}l(t)$\textit{ }$(t\rightarrow\infty)$.

(1) Assume $\omega_{q}<-1$. From $\lim\sup_{t\rightarrow\infty}$ $\log
q(t)/\log t<-1$, we obtain $\int_{0}^{\infty}q(t)dt$ $<\infty$, $\int_{0}%
^{x}q(t)dt\in R_{0}$, and $\omega_{S}=0$. Next, assume $\omega_{q}=-1$. If
$\omega_{S}\neq0,$ we have $\omega_{q}=\omega_{S}-1=-1$, which is a contradiction.

(2) Assume $\omega_{q}>-1$ and $\omega_{S}\neq0$. Then, we find that
$\omega_{q}=\omega_{S}-1<0$ and $\lim_{m\rightarrow E^{-}}\underline{u}_{\max
}/\underline{E}$ $=-\omega_{q}/(1-\omega_{q})$. Next, assume $\omega_{q}>-1$
and $\omega_{S}=0$. Then, $q(t)=o(t^{-1}l(t))$\textit{ }$(t\rightarrow\infty)$
and $\omega_{q}\leq-1$, thus contradicting the assumption.\hfill$\square$

\bigskip

\begin{center}
\textbf{References}\smallskip
\end{center}

\noindent\lbrack1] N. H. Bingham, C. M. Goldie and J. L. Teugels,
\textit{Regular variation} (Cambridge

University Press, Cambridge, 1987).

\noindent\lbrack2] N. H. Bingham and A. Inoue, `Ratio Mercerian theorems with
applications to

Hankel and Fourier transforms', \textit{Proc. London Math. Soc}. \textbf{79}
(1999), 626--648.

\noindent\lbrack3] N. H. Bingham and A. Inoue, `Tauberian and Mercerian
theorems for systems

of kernels', \textit{J. Math. Anal. Appl}. \textbf{252} (2000), 177--197.

\noindent\lbrack4] P. Cl\'{e}ment and J. Pr\"{u}ss, `Completely positive
measures and Feller semigroups',

\textit{Math. Ann.} \textbf{287} (1990), 73--105.

\noindent\lbrack5] W. Farkas, N. Jacob and R. L. Schilling, `Function spaces
related to continuous

negative definite functions: $\psi$-Bessel potential spaces',
\textit{Dissertationes Math}.

\textbf{393} (2001), 62 pp.

\noindent\lbrack6] W. Feller, \textit{An introduction to probability theory
and \ its application}, vol. I, II

(John Wiley and Sons, New York, 1957, 1966).

\noindent\lbrack7] F. Gesztesy and B. Simon, `Uniqueness theorems in inverse
spectral theory for

one-dimensional Schr\"{o}dinger operators', \textit{Trans. Amer. Math. Soc}.
\textbf{348} (1996),

349-373.

\noindent\lbrack8] F. Gesztesy and E. Tsekanovskii, `On matrix-valued Herglotz
functions', \textit{Math. }

\textit{Nachr}. \textbf{218} (2000), 61--138.

\noindent\lbrack9] Y. Hirashita, `Game pricing and double sequence of random variables',

Preprint, arXiv:math.OC/0703076 (2007).

\noindent\lbrack10] K. W. Homan, \textit{An analytic semigroup approach to
convolution Volterra equations}

(Delft University Press, Netherlands, 2003).

\noindent\lbrack11] N. Jacob and R. L. Schilling, `Subordination in the sense
of S. Bochner---an

approach through pseudo-differential operators', \textit{Math. Nachr}.
\textbf{178} (1996),

199--231.

\noindent\lbrack12] J. L. Kelly, `A new interpretation of information
rate',\ \textit{Bell System Tech. J.}

\textbf{35} (1956), 917--926.

\noindent\lbrack13] D. G. Luenberger, \textit{Investment science} (Oxford
University Press, Oxford, 1998).

\noindent\lbrack14] R. L. Schilling, `Subordination in the sense of Bochner
and a related functional

calculus', \textit{J. Austral. Math. Soc. Ser. A} \textbf{64} (1998), 368--396.

\noindent\lbrack15] D. V. Widder, \textit{The Laplace Transform }(Princeton
University Press, Princeton,

1941).

\bigskip

Chukyo University

Nagoya, 466-8666

Japan

e-mail: yukioh@cnc.chukyo-u.ac.jp
\end{document}